\documentclass[10pt]{article}

\usepackage{amsmath}
\usepackage{amssymb}
\usepackage{amsthm}
\usepackage{graphicx}
\usepackage{verbatim}
\usepackage{enumerate}

\font\elevenss=cmss11

\font\eightss=cmss8

\font\sixss=cmss8 at 6pt

\newfam\stfam
\textfont\stfam=\elevenss \scriptfont\stfam=\eightss
 \scriptscriptfont\stfam=\sixss
\def\st{\fam\stfam \elevenss}%
\theoremstyle{plain}
\newtheorem{thm}{Theorem}[section]
\newtheorem{lem}[thm]{Lemma}
\newtheorem{pr}[thm]{Proposition}
\newtheorem{cor}[thm]{Corollary}

\theoremstyle{definition}
\newtheorem{defns}[thm]{Definitions}

\newtheorem{example}[thm]{Example}

\newtheorem{alg}[thm]{Algorithm}
\theoremstyle{remark}
\newtheorem*{unremark}{Remark}
\newtheorem{remark}[thm]{Remark}
\newtheorem*{unremarks}{Remarks}
\newtheorem{remarks}[thm]{Remarks}
\numberwithin{equation}{section}

\newcommand{\Em}[1]{\textbf{#1}}

\def\grad{\nabla}
\def\ee{\epsilon}

\def\zz{{\bf z}}
\def\rr{{\bf r}}
\def\rhat{{\hat{\bf r}}}

\def\C{{\mathbb C}}
\def\Z{{\mathbb Z}}
\def\Q{{\mathbb Q}}
\def\N{{\mathbb N}}
\def\R{{\mathbb R}}

\def\Real{\; {\rm Re\,}}
\def\Imag{\; {\rm Im\,}}
\def\Arg{{\rm Arg\,}}

\def\sing{{\mathcal V}}
\def\M{{\mathcal M}}
\def\CC{{\mathcal C}}
\def\DD{{\mathcal D}}

\def\saddles{{\mbox{\st saddles}}}
\def\xns{W}
\def\bad{{\mbox{\st bad}}}
\def\degen{{\mbox{\st monkey}}}
\def\paths{{\mbox{\st paths}}}
\def\done{{\mbox{\st done}}}
\def\classify{{\mbox{\st subclass}}}
\def\class{{\mbox{\st class}}}

\def\one{{\bf 1}}
\def\nbd{{\mathcal N}}

\def\Res{{\rm Res }\,}

\def\singg{\sing'}
\def\Cdd{(\C^*)^d}

\def\At{{\tilde{A}}}
\def\critval{c_*}
\def\critalt{c_{\rm xy}}
\def\singx{X}
\def\singy{Y}
\def\rind{{\cal W}}
\def\cmax{c_{\max}}
\def\cmin{c_{\min}}
\def\Cox{\hfill \Box}
\def\ball{{\cal B}}
\def\algeb{{\mathbb A}}
\def\float{{\mathbb F}}

\makeatletter \def\romenumi{ \def\theenumi{\roman{enumi}}
\def\p@enumi{\theenumi} \def\labelenumi{(\@roman\c@enumi)}} 
\makeatother

\begin{document}

\begin{titlepage}

\begin{center}
{\large \bf Automatic asymptotics for coefficients of smooth, 
bivariate rational functions } \\
\end{center}
\vspace{5ex}
\begin{flushright}
Timothy DeVries\footnote{Department of Mathematics, University of 
Pennsylvania, 209 South 33rd Street, Philadelphia, PA 19104, 
{\tt tdevries@math.upenn.edu}}
  ~\\
Joris van der Hoeven\footnote{CNRS, Laboratoire LIX, \'Ecole Polytechnique, 
F-91228 Palaiseau Cedex, France,
{\tt vdhoeven@lix.polytechnique.fr}}$^,$\footnote{Research supported 
in part by the ANR-09-JCJC-0098-01 \textsc{MaGiX} project,
as well as Digiteo 2009-36HD grant and R\'egion Ile-de-France}
  ~\\
Robin Pemantle\footnote{Department of Mathematics, University of 
Pennsylvania, 209 South 33rd Street, Philadelphia, PA 19104, 
{\tt pemantle@math.upenn.edu}}$^,$\footnote{Research supported 
in part by National Science Foundation grant \# DMS 0905937}
\end{flushright}

\begin{abstract}
We consider a bivariate rational generating function
$$F(x,y) = \frac{P(x,y)}{Q(x,y)} = \sum_{r,s \geq 0} a_{rs} x^r y^s$$
under the assumption that the complex algebraic curve $\sing$
on which $Q$ vanishes is smooth.  Formulae for the asymptotics 
of the coefficients $\{ a_{rs} \}$ are derived in~\cite{PW1}.  
These formulae are in terms of algebraic and topological invariants
of $\sing$, but up to now these invariants could be computed only
under a minimality hypothesis, namely that the dominant saddle 
must lie on the boundary of the domain of convergence.  
In the present paper, we give an effective method for computing 
the topological invariants, and hence the asymptotics of 
$\{ a_{rs} \}$, without the minimality assumption.  This 
leads to a theoretically rigorous algorithm, whose 
implementation is in progress at
{\tt http://www.mathemagix.org}
\end{abstract}

\vfill

\noindent{Subject class: 05A15}

\end{titlepage}

\section{Introduction}

Consider a power series $F(\zz) = \sum_\rr a_\rr \zz^\rr$,
where $\rr$ varies over integer vectors in the orthant $(\Z^+)^d$ 
and $\zz^\rr$ denotes the monomial $z_1^{r_1} \cdots z_d^{r_d}$.  
We say that $F$ is the (ordinary) generating function for the 
array $\{ a_\rr \}$.  In analytic combinatorics, our aim is to 
derive estimates for $a_\rr$, given a simple description of $F$ 
as an analytic function.  An apparatus for doing this for various
classes of function is developed in a series of 
papers~\cite{PW1,PW2,BP-cones,pemantle-mvGF-AMS}; see also the
survey~\cite{PW9}.  Most known results on multivariate asymptotics
concern either rational functions or quasi-powers.  In the case
of rational functions $F(\zz) = P(\zz) / Q(\zz)$, the analysis 
centers on geometric properties of the pole variety 
$\sing := \{ \zz : Q(\zz) = 0 \}$.  

When $\sing$ has singularities, tools are required such as 
iterated residues, resolution of singularities or generalized 
Fourier transforms, and analyses exist only for specific examples.  
When $\sing$ is smooth, 
if certain degeneracies are avoided, asymptotic formulae for
$a_\rr$ may be given in terms of certain algebraic and topological
invariants of $\sing$.  Denote normalized vectors by
$$\rhat := \frac{\rr}{|\rr|} = \left ( \frac{r_1}{\sum_{j=1}^d r_j} , 
   \ldots , \frac{r_d}{\sum_{j=1}^d r_j} \right ) 
   := (\hat{r}_1 , \ldots , \hat{r}_d) \, .$$
Asymptotic formulae for $a_\rr$ depend on the direction $\rhat$
in which $\rr$ is going to infinity.  For example, Theorem~3.9 
of~\cite[Theorem~3.9]{pemantle-mvGF-AMS} gives an asymptotic 
formula for $a_\rr$ in terms of a sum over a set $\Xi$
of quantities that are easily computed via standard saddle
point techniques.  The set $\Xi = \Xi (\rhat)$ is a subset 
of the set $\saddles$ of saddle points of the function 
$h: \sing \to \R$ defined by
\begin{equation} \label{eq:h}
h (\zz) := h_{\rhat} (\zz) := - \sum_{j=1}^d \hat{r}_j \log |z_j| \, .
\end{equation}
The set $\saddles$ is readily computed but membership in the 
subset $\Xi$ is not easily determined, and in fact there is 
no known algorithm for doing so.  

The main result in this paper is a characterization of $\Xi$
when $d=2$ and the complex algebraic curve $\sing$ is smooth.
This leads to a completely effective and rigorous algorithm
for asymptotically computing $a_\rr$. Precise statements of
these results are somewhat technical, so we refer to the main text.
The organization of the paper is as follows.

Section~\ref{sec:background} reviews
background results from elsewhere which reduce the
computation of asymptotic formulae to identification of the set
$\Xi$ and of some path segments through each $\sigma \in \Xi$.
In particular, Section~\ref{sec:background} begins with the
integral representation of the general coefficient $a_\rr$,
Section~\ref{ss:residue} defines the residue of a meromorphic
form, Section~\ref{ss:intersection} reduces the integral for
$a_\rr$ to a lower-dimensional integral with some parameters
yet to be specified, and Section~\ref{ss:morse} selects a chain
of integration for this integral that results in an explicit
asymptotic estimate in terms of $\saddles$, $\Xi$ and
the so-called critical height $\critval$ (Theorem~\ref{th:ineffective}).

The new material begins in Section~\ref{sec:ident}.  Beginning
with the topology of $\sing$ near the coordinate axes
(Theorem~\ref{th:classify}), we give a topological characterization
of the critical minimax height $\critval$ to which the cycle of 
integration can be lowered (Theorem~\ref{th:critval}) and of the 
set $\Xi$ of contributing saddles (Theorem~\ref{th:Xi}).
Section~\ref{sec:computation} demonstrates how these  topological
computations of $\critval$ and $\Xi$ can be made completely
effective (Algorithm~\ref{alg:1} and Theorem~\ref{th:path}).
Details of implementation are given, as well as a discussion
of uniformity and boundary cases.

\section{Background} \label{sec:background}

In this section we review the framework for deriving estimates
for $a_\rr$, beginning with results for general $F$, then 
specializing to rational functions, and smooth pole variety 
$\sing$.  \Em{This section may be skipped by readers who want only
to understand the effective computation and not the underlying
analysis.}  Most of the material in this section is valid for any 
number of variables.  When we need to assume a bivariate function, 
we will state this but will also change the notation to use 
$x$ and $y$ in place of $z_1$ and $z_2$ and $(r,s)$ in place 
of $\rr$.

Computation of $a_\rr$ via complex analytic methods begins
with the multivariate Cauchy integral formula.  
\begin{equation} \label{eq:cauchy 1}
a_\rr = \frac{1}{(2\pi i)^{d}} \, \int_T
   \frac{F(\zz)}{z_1 \cdot \ldots \cdot z_d}
   \zz^{- \rr} \, d\zz  \, . 
\end{equation}
Here the torus $T$ is any product of a circles in each 
coordinate sufficiently small so that the product of the 
corresponding disks lies completely with the domain of 
holomorphy $\DD$ of $F$.  This version of Cauchy's formula 
can be found in most textbooks presenting complex analysis 
in a multivariable setting, and follows easily as an iterated 
form of the single variable formula; see, for 
example,~\cite[page~19]{shabat}.
The integrand is may be written
as $\zz^{-\rr} \omega$ where
\begin{equation} \label{eq:omega}
\omega := \omega_F := 
   \frac{F(\zz)}{z_1 \cdot \ldots \cdot z_d} \, d\zz  \, .
\end{equation}

Let $\Cdd := (\C \setminus 0)^d$ denote the set of complex
$d$-vectors with all nonzero coordinates and denote $\singg :=
\sing \cap \Cdd$.  A key step in obtaining asymptotic estimates
for $a_\rr$ is to transform the Cauchy integral equation
$(2 \pi i)^d a_\rr = \int_T \zz^{-\rr} \omega_F$ via an 
identity valid for any meromorphic form with a 
simple pole:
\begin{equation} \label{eq:residue}
\int_T \zz^{-\rr} \omega = 2 \pi i \, \int_\alpha \Res (\zz^{-\rr} \omega) 
\end{equation}
where $\Res (\cdot)$ is the residue operator defined in 
Section~\ref{ss:residue} and $\alpha$ is the intersection 
class defined in Section~\ref{ss:intersection}.  Putting
these together and specializing to $d=2$ leads to 
Lemma~\ref{lem:transfer} below:
$$a_\rr = \frac{1}{2 \pi i} 
   \, \int_{\alpha} x^{-r} y^{-s} \, \Res (\omega) \, .$$

\subsection{The residue form} \label{ss:residue}

The residue form $\Res (\eta)$ is a holomorphic $(d-1)$-form 
on $\singg$.  The specification of this form will not be important 
for the subsequent analysis, but for completeness we include
the following definition.  If $\eta = (P/Q) \, d\zz$ is any 
meromorphic form on a domain $U$ with simple pole on
a set $\sing \subseteq U$, then we define 
\begin{equation} \label{eq:res 1}
\Res(\eta) := \iota^* \theta \, 
\end{equation}
where $\theta$ is any solution to
\begin{equation} \label{eq:res 2}
dQ \wedge \theta = P \, d\zz \, .
\end{equation}
Existence and uniqueness are well known, and a proof
can be found in~\cite[Proposition~2.6]{devries09} among
other places (see, e.g.,~\cite{aizenberg-yuzhakov}).
One important property of the residue form is that if
$g: \sing \to \C$ is holomorphic then
$$\Res (g \cdot \eta) = g \cdot \Res(\eta) \, .$$
In particular, $\Res (\zz^{-\rr} \, \omega ) = \zz^{-\rr} 
\Res (\omega)$.

\subsection{The intersection class} \label{ss:intersection}

Because $\Res (\zz^{-\rr} \, \omega)$ is holomorphic on $\singg$, 
the integral over a $(d-1)$-cycle $B \subseteq \singg$ depends 
only on the homology class of $B$ in $H_{d-1} (\singg)$.
The chain of integration $\alpha$ in~\eqref{eq:residue} is
really a homology class known as the \Em{intersection class}
of $T$ with $\singg$.  The construction of the intersection
class is quite general and may be found in the literature.
Conditions are given in~\cite[Appendix~A]{PW-book} for the
intersection class to be uniquely defined.  When it is not,
however, the integral over $\alpha$ must be the same for
any choice of intersection class, $\alpha$, so we will not
pursue it further here.

Because we will need an explicit construction of a cycle in this
class, we will give a quick construction of $\alpha$ and derivation 
of~\eqref{eq:residue}.  
We begin with a form of the Cauchy-Leray 
residue theorem, which may be found in~\cite[Theorem~2.8]{devries09}.
\begin{lem}[Cauchy-Leray Residue Theorem] \label{lem:cauchy-leray}
Let $\eta$ be a meromorphic $d$-form on domain $U \subseteq \C^d$ 
with pole variety $\sing \subseteq U$ along which $\eta$ has only 
simple poles.  Let $N$ be a $d$-chain in $U$, locally the product 
of a $(d-1)$-chain $C$ on $\sing$ with a circle $\gamma$ in the 
normal slice to $\sing$, oriented as the boundary of a disk 
oriented positively with respect to the complex structure of 
the normal slice.  Then 
$$\int_N \eta = 2\pi i \int_C \Res(\eta) \, .$$
$\Cox$
\end{lem}

\begin{cor} \label{cor:homot}
Let $F, P, Q$ and $\sing$ be as above and let $\eta$ be meromorphic.
Let $\M_1$ and $\M_2$ be compact $d$-manifolds in $\C^d$ and let $H$ 
be a smooth homotopy from $\M_1$ to $\M_2$ that intersects $\sing$ 
transversely in a compact manifold $\CC$.  Then 
$$\int_{\M_1} \eta - \int_{\M_2} \eta 
   = 2 \pi i \, \int_{\CC} \Res (\eta) \, .$$
$\Cox$
\end{cor}

\noindent{\sc Proof:} Let $N$ be the boundary of a tubular neighborhood
$W$ of $\CC$ in $\C^d$.  Then $\eta$ is holomorphic on $H \setminus W$.
The boundary of $H \setminus W$ is $\M_2 - \M_1 + N$, and $d \eta$ 
vanishes on $H \setminus W$ because the differential of any holomorphic
$d$-form vanishes on $\C^d$.  Tubular neighborhoods of a smooth algebraic
hypersurface always have a product structure (see, 
e.g.,~\cite[Proposition~A.4.1]{PW-book}).  Therefore, Stokes' Theorem 
and Lemma~\ref{lem:cauchy-leray} together yield
$$\int_{\M_1} \eta = \int_{\M_2} \eta + \int_N \eta
   = \int_{\M_2} \eta + 2 \pi i \, \int_{\CC} \Res (\eta) \, .$$
$\Cox$

\begin{remarks} \label{rem:orient}
One sees that the orientation of $\CC$ must be chosen so that 
the boundary of its product with a positively oriented disk 
in the normal slice is homologous to $\M_1 - \M_2$.  Also
note that result is true when $H \setminus W$ is any 
cobordism: we used only that $\partial (H \setminus W) = 
\M_2 - \M_1 + N$, not that $H$ was a homotopy. 
\end{remarks}

Let $T_\ee$ be the torus $\{ |z_1| = \cdots = |z_d| = \ee \}$ and let
$T_{\ee , L}$ be the torus $\{ |z_1| = \cdots = |z_{d-1}| = \ee \, ,
|z_d| = L \}$.  Fix $\ee > 0$ and suppose that for sufficiently 
large $L$, the torus $T_{\ee , L}$ does not intersect $\sing$.
We claim that for $r_d$ sufficiently large,
\begin{equation} \label{eq:vanish}
\int_{T_{\ee , L}} \zz^{-\rr} \, \omega = 0 \, .
\end{equation}
Indeed, 
\begin{eqnarray*}
|\zz^{-\rr} \omega| & \leq & \ee^{-\sum_{j=1}^{d-1} r_j} L^{-r_d}
   \sup_{\zz \in T_{\ee,L}} |\omega (\zz)| \\[1ex]
& \leq & C L^{-r_d} L^\beta
\end{eqnarray*}
where $\beta$ is the degree of the rational form $\omega$.
The chain of integration has volume $O(L)$, from which we 
conclude that $\int_{T_{\ee , L}} \zz^{-\rr} \, \omega = 
O(L^{1 + \beta - r_d})$.  Once $r_d > 1 + \beta$, we therefore
have $\int_{T_{\ee , L}} \zz^{-\rr} \, \omega = o(1)$ as
$L \to \infty$.  But $T_{\ee , L}$ avoids $\sing$ for sufficiently
large $L$, hence by Stokes' Theorem, $\int_{T_\ee , L} \zz^{-\rr} 
\, \omega$ is constant for sufficiently large $L$, hence zero,
establishing~\eqref{eq:vanish}.

The intersection class may now be constructed as follows.  
Let $H_L$ be the homotopy defined by 
$$H_L (\zz , t) = (z_1 , \ldots , z_{d-1} , 
   (1 - t) z_d + t (L/\ee - 1) z_d) \, .$$
Suppose that $H_L [T_\ee]$ intersects $\sing$ transversely.  We let
\begin{equation} \label{eq:alpha}
\alpha = [\CC] ~~~ \mbox{ where } ~~~ \CC := H_L [T_\ee] \cap \sing \, .
\end{equation} 
This is independent of $L$ once $L$ is sufficiently large.  
We take the orientation of $\CC$ as in the remark following
Corollary~\ref{cor:homot}.

Summing up the residue and intersection cycle constructions
we have:
\begin{lem}[Residue representation of $a_\rr$] 
\label{lem:transfer}
Suppose $d=2$.  Then for sufficiently small $\ee > 0$,
\begin{equation} \label{eq:transfer}
a_\rr = \frac{1}{(2 \pi i)^{d-1}} 
   \; \int_{\alpha} x^{-r} y^{-s} \, \Res (\omega) 
\end{equation}
holds with $\Res$ defined as in~\eqref{eq:res 1}~--~\eqref{eq:res 2} 
and the intersection class $\alpha$ defined in~\eqref{eq:alpha}.
\end{lem}

\noindent{\sc Proof:} 
Let $S$ be the set of $x \in \C$ such that there is a 
sequence $(x_n , y_n)$ in $\sing$ with $x_n \to x$ and
$y_n \to \infty$.  The set $S$ is finite and therefore
$\min \{ |x| : x \in S , x \neq 0 \}$ is strictly positive.
Taking $\ee$ less than this minimum guarantees that
for sufficiently large $L$, the torus $T_{\ee , L}$ does
not intersect $\sing$.  The set of $(x,y) \in \sing$
where $\partial Q / \partial y$ vanishes is also finite.
Again, taking $\ee$ less than the least nonzero modulus
of such a point guarantees that the homotopy $H_L$ 
intersects $\sing$ transversely.  Having verified these
two suppositions, the construction of $\alpha$ is completed.
Now apply Corollary~\ref{cor:homot} with $\M_1 = T_\ee,
\M_2 = T_{\ee , L}, N = \CC$ and $\eta = \zz^{-\rr} \omega$,
where $\ee$ is small enough so that $T_\ee \subseteq \DD$.
Using~\eqref{eq:vanish} to see that $\int_{\M_2} \eta = 0$ 
and combining with Cauchy's integral formula gives
$$(2 \pi i)^d a_\rr = \int_{T_\ee} \zz^{-\rr} \omega 
   = 2 \pi i \, \int_{\alpha} \Res (\zz^{-\rr} \omega) 
   = 2 \pi i \, \int_{\alpha} \zz^{-\rr} \, \Res (\omega) \, .$$
$\Cox$

\begin{remark} \label{rem:d=2}
The only use thus far of the restriction to $d=2$ is in verifying
that $H_L$ intersects $\sing$ transversely and $T_{\ee , L}$ not
at all.  Transversality may be accomplished for any $d$ via a
perturbation and it is not hard to handle a nonempty intersection
of $\sing$ with $T_{\ee , L}$, so no essential use has been made
yet of the restriction to $d=2$.
\end{remark}

\subsection{Morse theory} \label{ss:morse}

Saddle point integration techniques require that we choose a
chain of integration $\CC'$ in the class $\alpha$ to (roughly)
minimize the modulus of the integrand $\zz^{-\rr} \omega$.
It will suffice to minimize the factor of $\zz^{-\rr}$,
the form $\omega$ being bounded on compact chains.  The log
modulus of $\zz^{-\rr}$ is $\sum_{j=1}^d r_j \log |z_j|$.
Letting $n = |\rr| = \sum_{j=1}^d r_j$, we have
$$\log \left | \zz^{-\rr} \right | = 
   \sum_{j=1}^d r_j \log |z_j| =
   n h_{\rhat} (\zz)$$
which explains our definition of $h$ in~\eqref{eq:h}.  It will
be helpful to consider the manifold $\singg$ from a Morse
theoretic viewpoint with height function $h$.  

\begin{defns}
Let $\M$ be any manifold and let $h: \M \to \R$ be any smooth, 
proper function with isolated critical points (a critical point 
is where $\grad h$ vanishes).  Let $\sing^{\geq c}$ denote the 
subset of $\sing$ consisting of points $\zz$ with $h(\zz) \geq c$.  
Define $\sing^{> c}, \sing^{\leq c} , \sing^{< c}$ and $\sing^{= c}$ 
similarly.  Define the set $\saddles = \saddles (\M , h)$ to be the 
set of critical points for $h$.  The set $h[\saddles]$ is called
the set of critical values.  For fixed $(r ,s)$, the set $\saddles$ 
is a zero-dimensional variety, whose defining equations we call the 
\Em{critical point equations}:
\begin{eqnarray}
Q(x,y) & = & 0 \label{eq:crit} \\
s \, x \, Q_x - r \, y \, Q_y & = & 0 \, . \nonumber
\end{eqnarray}
Here and henceforth, we use subscript notation for partial derivatives.
\end{defns}
When $h : \singg \to \R$ is the height function $h_{\rhat}$ 
above, by convention we extend $h$ to $\sing$ taking
$h(0,y) = h(x,0) = + \infty$, so that points of 
$\sing \setminus \singg$ are in every $\sing^{\geq c}$.

\begin{thm} \label{th:morse}
Let $h : \M \to \R$ be a smooth, proper function with isolated
critical points.  
\begin{enumerate} \romenumi
\item Suppose there are no critical values of $h$ in the interval
$[a,b]$.  Then the gradient flow induces a homotopy in $\M$ between
$\M^{\leq b}$ and $\M^{\leq a}$, and also between $\M^{\geq b}$
and $\M^{\geq a}$.  This homotopy carries $\M^{= b} = \partial
\M^{\geq b} = \partial^{\M \leq b}$ to $\M^{= a}$.
\item Let $a < c < b$ such that $c$ is the unique critical value
of $h$ in $[a,b]$.  Let $\sigma_1 , \ldots , \sigma_k$ enumerate
the critical points of $h$ at height $c$.  Then $H_* (\M^{\leq b} ,
\M^{\leq a})$ decomposes naturally as the direct sum of 
$H_* (\M_{\sigma_j}^{\leq b} , \M_{\sigma_j}^{\leq a})$, 
where $\M_{\sigma_j}$ is the intersection of $\M$ with a 
sufficiently small neighborhood of $\sigma_j$.
\item Let $\CC$ be a cycle in $\M$ and let $c_*$ be the infimum of
$c$ such that $\CC$ is homologous in $\M$ to a cycle in $\M^{\leq c}$.
Then $c$ is either $-\infty$ or a critical value, and in the latter
case, for sufficiently small $\ee > 0$, the cycle $\CC$ projects to 
a nonzero homology class in $H_* (\M , \M^{\leq c - \ee})$.
\end{enumerate}
\end{thm}

\noindent{\sc Proof:} Part~$(i)$ is the first Morse lemma,
a proof of which may be found for example in~\cite[Theorem~3.1]{milnor}
or~\cite{GM}.  A generic perturbation separates the critical values, 
one may then prove part~$(ii)$ by a standard excision argument. 
Part~$(iii)$ is proved in~\cite[Lemma~8]{pemantle-mvGF-AMS}.
$\Cox$

\begin{remark}
Classical Morse theory combines the first Morse lemma (no 
topological change between critical values) with a description
of the change at a critical value (an attachment map).  The
results we have quoted require only the first Morse lemma. 
For this it is not necessary to assume $h$ is a Morse function,
but only that it is a smooth, proper map with isolated critical 
points.  
\end{remark}

\begin{example}[bivariate harmonic case] \label{eg:harmonic}
Suppose $d=2$.  Let $\sigma$ be a critical point of $h$ on $\sing$ 
at height $c$ and suppose that locally $h$ is the real part of an 
analytic function.  Then there is an integer $k \geq 2$ and there
are local coordinates $\psi: U \to \sing$ on a neighborhood of 
the origin in $\C^1$ such that $\psi (0) = \sigma$ and
$$h (\psi (z)) = c + \Real \{ z^k \} \, .$$
Figure~\ref{fig:k-saddle} illustrates this for $k=4$.  The shaded 
regions are where $h > c$.  On the right, the connected local
components of $\{ h > c \}$ are number counterclockwise from~1
to~$k$ and the local components of $\{ h < c \}$ are numbered
correspondingly.
\begin{figure}[ht]
\centering
\includegraphics[scale=0.40]{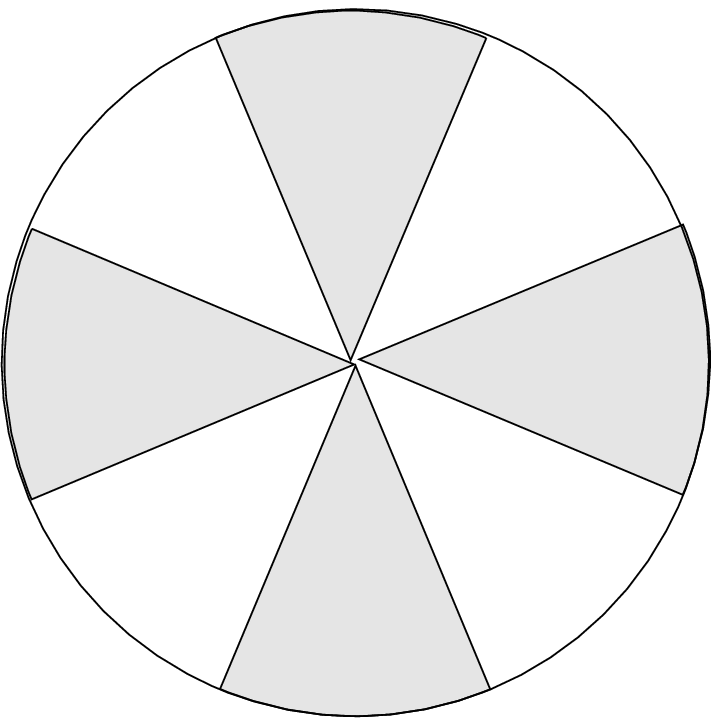} \hspace{1.5in}
\includegraphics[scale=0.40]{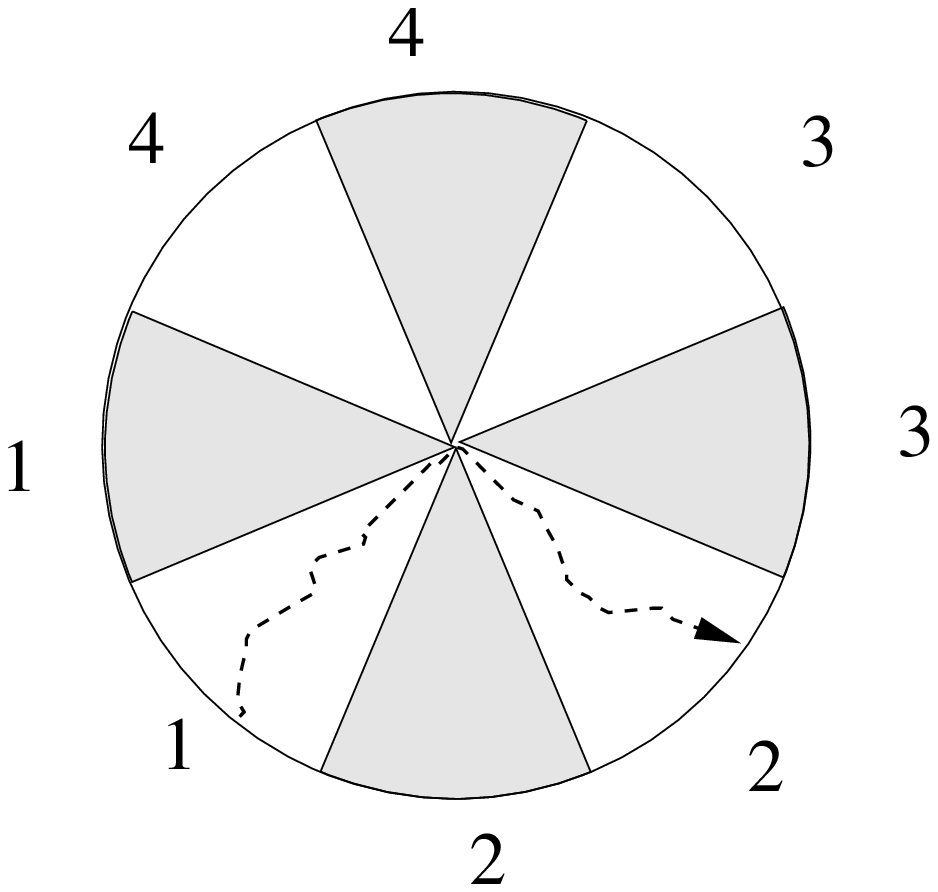}
\caption{A saddle of order~4 and the path $\gamma_{1,2}$}
\label{fig:k-saddle}
\end{figure}
The pair $\tilde{\M} := (\sing^{\leq c + \ee} \cap \psi[U] , 
\sing^{\leq c + \ee} \cap \psi[U])$ is homotopy equivalent to 
a disk modulo $k$ boundary arcs.  Let $\gamma_j$ denote an arc 
from the boundary to the center that remains in the $j^{th}$ 
local component of $\{ h < c \}$.  Then 
$H_1 (\tilde{\M})$ is generated by the relative 1-cycles
$\gamma_i - \gamma_j$.  We have the relations
$\sum_{i,j} b_{i,j} (\gamma_i - \gamma_j) = 0$ if and only
$$\sum_j b_{i,j} = \sum_j b_{j,i}$$ 
for every $i \leq k$.  In particular, $\{ \gamma_i - \gamma_{i+1} \; : \;
1 \leq i \leq k-1 \}$ is a basis for $H_1 (\tilde{\M})$.

We remark that our assumptions imply that $h$ has no local extrema 
on $\M$.  Therefore, as $c$ decreases, the connected components 
of $\M^{\geq c}$ can merge but new components never arise.
\end{example}
This leads to a canonical form in dimension $d=2$ for classes in
$H_1 (\singg)$, choosing a representing cycle so as to be in
position for saddle point integration.
\begin{pr} \label{pr:d=2}
Let $h = h_{\rhat}$.  Let $\gamma$ be any cycle in $H_1 (\singg)$.
Then either $\gamma$ is homologous to a cycle supported at arbitrarily
small height or there is a critical value $c_*$ and a nonempty set 
$\xi \subset \saddles$ of critical points at height $c_*$ 
such that the following holds.
\begin{quote}
The cycle $\gamma$ is homologous in $H_1 (\singg)$ to a cycle $\CC'$
in $\sing^{\leq c_*}$; the maximum of $h$ on $\sing$ is achieved
exactly on $\Xi$; in a neighborhood of each $\sigma \in \Xi$ of
order $k = k(\sigma) \geq 2$, there is a subset $\paths (\sigma)$ of
$\{ 1 , \ldots , k (\sigma) - 1 \}$ such that the restriction of
$\CC'$ to the neighborhood is equal to $\sum_{i \in T (\sigma)}
\gamma_i - \gamma_{i+1}$.
\end{quote}
\end{pr}

As promised, integrals of the residue form $\Res (x^{-r} y^{-s} 
\omega_F)$ are known.  To summarize this, we let $F = P/Q$ be 
rational, with $\omega$ defined in~\eqref{eq:omega} and
$\Res (\omega)$ defined in Section~\ref{ss:residue}.  With
$\gamma_i$ as in Figure~\ref{fig:k-saddle}, we then have the
following explicit integration estimates.  Note that the 
magnitude of $x^{-r} y^{-s}$ is $\exp [ (r+s) \critval ]$
when $(x,y) \in \Xi$, which makes the remainder term 
exponentially smaller than the summands.
\begin{lem}[saddle point integrals]
\label{lem:integrals}
Let $\sigma = (x_0 , y_0) \in \saddles$, satisfy the critical
point \hfill \\ equations~\eqref{eq:crit}.  Denote the order of 
vanishing of $h$ at $(x_0 , y_0)$ by $k = k(x_0 , y_0)$.
Then \hfill \\ $(2 \pi i ) \int_{\gamma_i - \gamma_{i+1}} 
\Res (x^{-r} y^{-s} \omega)$ is given by the formula
$\Phi (\sigma , i)$ where
\begin{equation}  \label{eq:integrals}
\Phi (\sigma , i) := 
   \frac{K(\sigma , i)}{\sqrt{2 \pi (r+s)}}\,  x_0^{-r} y_0^{-s} 
   + O \left ( \exp [(\critval - \ee) (r+s)] \right )  
\end{equation} 
where $K(\sigma , i)$ is a constant term given by~\eqref{eq:K}
when $d=2$ and more generally by an explicit formula which may
be written in a number of ways, e.g.~\cite[Theorem~3.5]{PW1}
or~\cite[Theorem~3.3]{BBBP}.
When $k=2$, the remainder term is uniform as $(x_0 , y_0 , r , s)$
varies with $(x_0 , y_0) \in \saddles (r,s)$ and $(x_0 , y_0)$
varying over a compact set on which $k$ is constant.  
When $k > 2$, the estimate is uniform for fixed $(x_0 , y_0)$ 
and $(r,s) \to \infty$ with $r = \lambda s + O(1)$, 
where $\lambda = x Q_x / (y Q_y)$ evaluated at $(x_0 , y_0)$.  
$\Cox$
\end{lem}

\begin{remarks} \label{rem:airy} ~~\\[-3ex]
\begin{enumerate} \romenumi
\item In fact a uniform estimate holds for $r = \lambda s + O(s^{1/2})$ 
although this is more complicated and involves Airy functions.
\item When $d=2$ the constant $K(\sigma , i)$ is given 
in~\cite[Theorem~3.1]{PW1} as
\begin{equation} \label{eq:K}
K(\sigma , i) = P(x , y) \sqrt{ \frac{r+s}{s} \, 
   \frac{- Q_y}{x \, \psi } }
\end{equation}
evaluated at $(x_0 , y_0)$, 
where $\psi (x,y)$ is the expression in equation~\eqref{eq:psi} below.
\item When $k > 2$, the constant $K(\sigma  , i)$ is the leading term
in a series expansion for $\lambda$ near $\sigma$.  The explicit
formula involves partial derivatives of $P$ and $Q$ up to order~$k$
and may be derived from, e.g.,~\cite[Theorem~3.3]{PW1}.
\end{enumerate}
\end{remarks}

Putting together Lemma~\ref{lem:transfer} with Proposition~\ref{pr:d=2}
and Lemma~\ref{lem:integrals} leads to the following asymptotic 
expression for $a_\rr$, which has appeared several times before;
for example, when $\sigma$ is a minimal point, this goes back 
to~\cite[Theorem~3.1]{PW1}; for an alternative formulation
see~\cite[Theorem~3.3]{BBBP}.
\begin{thm}[smooth point asymptotics] \label{th:ineffective}
$$ a_{rs} = \left ( \frac{1}{2 \pi i} \right )^{d-1} 
   \sum_{\sigma \in \Xi} \;\;\;  \sum_{i \in \paths (\sigma)}
   \Phi (\sigma , i) 
   + O \left ( \exp [(\critval - \ee) (r+s)] \right )  
   \, .$$
$\Cox$
\end{thm} \vspace{-0.2in}

This completes the summary of known results.  We now turn to 
the main work which is to identify and compute the quantities 
$\critval, \Xi$ and $\paths (\sigma)$ appearing in 
Theorem~\ref{th:ineffective} and formula~\eqref{eq:integrals}.

\section{Identifying $\critval$ and $\Xi$} \label{sec:ident}

We begin by constructing the intersection class.  We will
take a homotopy which keeps the torus small in the $x$-coordinate.
If $\sing$ intersects the $y$-axis transversely at $(0,y)$, then
the set $\sing \cap \{ |x| \leq \ee \}$ will have a 
component that is a topological disk nearly coinciding
with a component of $\sing^{\geq c}$ with $c = - \hat{r} 
\log \ee$.  To see that something like this holds more generally, 
we recall some facts about Puiseux expansions.  Consider the 
multivalued function $x \mapsto y(x)$ that solves $Q(x,y(x)) = 0$.  
When $x$ is in a sufficiently small punctured neighborhood of~0, 
there is a finite collection of convergent series of the form
\begin{equation} \label{eq:puiseux}
y(x) = \sum_{j \geq j_0} c_j x^{j/k} 
\end{equation}
where $j_0 \in \Z$, $k \in \Z^+$, and $c_{j_0} \neq 0$; such 
a series yields $k$ different values of $y$ for each $x$; 
together, the collection of series yields each solution 
to $Q(x,y) = 0$ exactly once.  We always assume that $k$ is
chosen as small as possible.  A proof of these facts may 
be found in~\cite[Theorem~VII.7]{flajolet-sedgewick-anacomb}. 
Our assumption that $Q(0,0) \neq 0$ implies $j_0 \leq 0$
in every Puiseux series solution.  Those with $j_0 = 0$
correspond to points at which $\sing$ intersects
the $x$-axis, the $y$-value being equal to the coefficient $c_0$.
Under the assumption that $\sing$ is smooth, for each such
$y$-value there will be only one series\footnote{In
fact the analyses in this paper are valid when self-intersections
of $\sing$ are allowed at points $(0,y)$, with the result
that the closures of the components of $\singg \cap \{ |x| \leq \ee \}$
may intersect.}.  If $k = 1$, the
intersection with the $y$-axis is transverse but it is possible
that $k \geq 2$ and $\sing$ is smooth but with non-transverse
intersection.  Finally, series with $j_0 < 0$ correspond to
components of $\sing \cap \{ |x| \leq \ee \}$ with $y$ going
to infinity as $x$ goes to zero.  We begin by showing that 
each series gives a single component and describing the
embedding of the component.  First, to clarify, we give an
example in which all these possibilities occur.

\begin{example} \label{eg:all}
Let $Q(x,y) = 1 - 3y + 2 y^2 - 6 x y^4 + x^3 y^5$.  It is easy 
to check that $\sing$ is smooth.  The total number of solutions
should be the $y$-degree of $Q$, namely~5.  Setting $x=0$ we 
find two points $(0,1)$ and $(0,1/2)$ of $\sing$ on the $y$-axis.  
These correspond to two series solutions with $j_0 = 0$ and $k=1$:
\begin{eqnarray*}
y(x) & = & 1 + 6 x + 72 x^2 + \cdots \\
y(x) & = & \frac{1}{2} - \frac{3}{8} x + \frac{45}{32} x^2 + \cdots 
\end{eqnarray*}
There are also two series solutions with $j_0 < 0$, one with 
$k=1$ and one with $k=2$:
\begin{eqnarray*}
y(x) & = & \frac{1}{\sqrt{3}} x^{-1/2} - \frac{3}{4}  
   + \frac{19 \sqrt{3}}{32} x^{1/2} + \cdots \\
y(x) & = & 6 x^{-2} - \frac{1}{18} x + \frac{1}{72} x^3 \cdots  \, .
\end{eqnarray*}
The former of these yields two solutions.  Although we do not
need it here, we remark that these five solutions can be 
enumerated by counting $y$-displacement along the Northwest 
boundary of the Newton polygon: going up from $(0,0)$ to $(0,2)$
counts the first two solutions, going from $(0,2)$ to $(1,4)$
yields the two solutions $y \sim \pm 1/\sqrt{3x}$ and going from
$(1,4)$ to $(3,5)$ yields the solution $y \sim 6 x^{-2}$.
\end{example}

Let $\singx^{\geq c}$ (respectively $\singy^{\geq c}$ denote 
the union of those components of $\sing^{\geq c}$ containing 
arbitrarily small values of $x$ (respectively $y$).  
Say that the \Em{direction} of a Puiseux series solution to 
$Q(x,y(x)) = 0$ is $- j_0 / k$.  The direction is always nonnegative
because $j_0$ is nonpositive due to $Q(0,0) \neq 0$.  Let 
$$\bad := \{ \beta > 0 \, : \, \beta \mbox{ is a direction of a series 
   solution to } Q(x,y(x)) = 0 \}$$
denote the set of positive directions.

\begin{thm}[$\singg$ near the axes] \label{th:classify}
Fix $\rhat$ with $r/s \notin \bad$.  Suppose that $h$ has 
isolated critical points.  Then when $\ee > 0$ is sufficiently 
small and $c$ is sufficiently large, the following hold.
\begin{enumerate} \romenumi
\item The graph of each Puiseux expansion of $y(x)$ over the
punctured disk $\{ 0 < |x| < \ee \}$ is a component of
$\singg \cap \{ |x| \leq \ee \}$.  
\item Each such component is a topological disk and its boundary
winds $k$ times around the $y$-axis. 
\item The components of $\singx^{\geq c}$ are exactly the graphs 
over some neighborhood of~0 of those Puiseux series whose directions 
are less than $r/s$.
\item Components of $\singg \cap \{ |x| < \ee \}$ and $\singx^{\geq c}$
with the same Puiseux series are homotopic in $\singg$.
\end{enumerate}
\end{thm}

\noindent{\sc Proof:}
Fix any Puiseux series solution $\sum_{j \geq j_0} c_j x^{j/k}$ 
on a punctured disk $U$ and define $g(x) := \sum c_j x^j$.  From
this we define the function 
\begin{align*}
G:U&\longrightarrow D\\
x&\longmapsto (x^k,g(x)) \, .
\end{align*}
We will show that $G$ is a diffeomorphism onto its image.  For 
sufficiently small $\ee$ the images of the punctured $\ee$-disk
under different Puiseux series are disjoint, which will finish
the proof of~$(i)$ and~$(ii)$.  

First we check that $G$ is one-to-one on a sufficiently small punctured
disk.  If not then we may find arbitrarily small $x_1$ and $x_2$
such that $x_1^k = x_2^k$ and $g(x_1) = g(x_2)$.  This means that
there is some $k^{th}$ root of unity $\xi$ and such that
$g(x) = g(\xi x)$ for arbitrarily small $x$.  For some $m$,
the functions $x^m g(x)$ and $x^m g(\xi x)$ are holomorphic,
and if they agree on a sequence converging to zero then they
must define the same function on $U$.  It follows that their 
Taylor expansions agree, and hence that the only nonzero
coefficients are those indexed by a multiple of $k$.  This
contradicts the minimality of $k$ in the Puiseux series
representation and we conclude that $g$ is one-to-one on a
sufficiently small disk.  For any smooth function $\phi$, 
the function $(x^k , \phi (x)) \mapsto x$ is diffeomorphism
near any $x \neq 0$ because $x^k \mapsto x$ is locally smooth
except at~0.  This proves that $G$ is a diffeomorphism onto
its image, establishing~$(i)$ and~$(ii)$.  To continue, we need
a lemma.

\begin{lem} \label{lem:radial}
When $\ee$ is sufficiently small, for any $\theta$ the radial arcs
$$t \mapsto (t e^{i \theta} , y(t e^{i \theta})) , \; 0 \leq t \leq \ee$$ 
are monotone in height, that is, 
$$h_\theta (t) := h (t e^{i \theta} , y(t e^{i \theta}))$$ 
is monotone increasing or decreasing.  Furthermore, as
$t \downarrow 0$, $h_\theta (t) \to \pm \infty$ according 
to whether the direction $\beta = -j_0 / k$ of the Puiseux 
series is less or greater than $r/s$.
\end{lem}

\noindent{\sc Proof:} Write
$$y(x) = c x^{-\beta} (1 + \phi (x))$$
where the remainder $\phi$ satisfies $\phi (x) = o(1)$
and $x \phi' (x) = o(1)$ as $x \to 0$.
Plugging this into the function $H(x,y) := - \hat{r} \log x 
- \hat{s} \log y$ gives
\begin{eqnarray}
\frac{d}{dx} H(x,y(x)) & = & \frac{\hat{s} \beta - \hat{r}}{x}
   - \frac{\phi' (x)}{1 + \phi (x)} \nonumber \\[1ex]
& = & \frac{1}{x} \left ( \hat{s} \beta - \hat{r} 
   - \frac{x \phi' (x)}{1 + \phi (x)} \right ) \nonumber \\[1ex]
& \sim & \frac{\hat{s} \beta - \hat{r}}{x} \label{eq:H}
\end{eqnarray}
because we have assumed $\hat{s} \beta - \hat{r} \neq 0$.
The function $h$ is the real part of $H$, which along with
the chain rule yields
$$\frac{d}{dt} h_\theta (t) = \cos \theta \Real \frac{d}{dx} H(x,y(x)) 
   - \sin \theta \Imag \frac{d}{dx} H(x,y(x))  \, .$$
Plugging in~\eqref{eq:H} then gives 
$$\frac{d \,  h_\theta (t)}{dt} \sim \frac{\hat{s} \beta - \hat{r}}{|x|}$$
which finishes the proof of the lemma.
$\Cox$

\noindent{\sc Proof of Theorem~\protect{\ref{th:classify}}, continued:}
Each component $B$ of $\singx^{\geq c}$ has, by definition, points
with arbitrarily small values of $x$ which are therefore on the
graph of one of the Puiseux series.  Taking a sufficiently small
neighborhood, the graphs of the series are disjoint, therefore
the map from component to series is well defined.  We need to
see that the range of the map from components to series is
exactly those with directions les than $r/s$.  If $B$
has direction less than $r/s$ then $h \to \infty$ as $x \to 0$
on the corresponding component, whence the corresponding 
component of $\singg \cap \{ |x| < \ee \}$ maps to $B$.
On the other hand, if the direction of $B$ is greater than
$r/s$ then $h \to -\infty$ as $x \to 0$, from which it follows 
that $\singx^{\geq c}$ cannot contain points in $B$ with
arbitrarily small $x$ values and hence such a component is not
the graph of a Puiseux series about $x=0$.  This proves~$(iii)$.

Let $B$ be a component of $\singx^{\geq c}$ where $c$ is sufficiently
large so that the projection onto $x$ is contained in a
sufficiently small disk so that Lemma~\ref{lem:radial} applies.
A radial homotopy in the $x$-coordinate (one that maps 
$x := \rho e^{i \theta}$ to $f(\theta) x$ and maps $y(x)$
to $y(f(\theta) x)$) then shrinks $B$ to the same component 
with any specified greater value of $c$, or the 
corresponding component of $\singg \cap \{ |x| \leq \ee \}$
with as long as $\ee$ is at most the in-radius of the
original $x$-neighborhood.  This establishes~$(iv)$ and
also that either type of component remains in the same 
homotopy class as $c$ or $\ee$ is varied. 
$\Cox$

Switching the roles of $x$ and $y$ in Theorem~\ref{th:classify},
it follows that for sufficiently large $c$, the set $\sing^{\geq c}$ 
is the disjoint union of $\singx^{\geq c}$ and $\singy^{\geq c}$.
Define $\critalt$ to the least value of $c$ for which this 
remains true, which is also that greatest value of $c$ such 
that $\sing^{\geq c}$ has a component containing both a point
$(0,y)$ and a point $(x,0)$.  By the first Morse lemma, if
$h$ is smooth and proper, the topology of $\sing^{\geq c}$ 
cannot change between critical values of $h$.  In particular 
$\critalt$ must be a critical value $c_j$.  The next two 
results complete the description of $\critval$ and $\Xi$.
\begin{thm}[characterization of $\critval$] \label{th:critval}
Assume that $\sing$ is smooth and that $h = h_{\rhat}$ is smooth 
and proper with isolated critical points.  Then \\[-4ex]
\begin{enumerate} \romenumi
\item For each $c \geq \critalt$, the cycle $\partial \singx^{\geq c}$
is in the homology class $\alpha$ from Lemma~\ref{lem:transfer}.
\item \hspace{0in} $\critval = \critalt$. 
\end{enumerate}
\end{thm}
For each critical point $\sigma$ with $h(\sigma) = \critval$ 
let $k = k (\sigma)$ denote its order, let $\nbd$ denote a
small neighborhood of $\sigma$ in $\sing$.  Label the components
of $\nbd \cap \{ h > \critval \}$ in counterclockwise order by
$R_1 , \ldots , R_k$ as in Figure~\ref{fig:k-saddle}.  Each such 
region $R$ will be called an $x$-region if $R \subseteq \singx^{> c}$ 
and a $y$-region if $R \subseteq \singy^{> c}$.  Because $c \geq \critalt$, 
these are mutually exclusive.
\begin{thm}[characterization of $\Xi$ and $\paths$] \label{th:Xi}
~~\\[-5ex]
\begin{enumerate} \romenumi
\item $\sigma \in \Xi$ if and only if the regions $R_1 , \ldots , R_k$
include at least one $x$-region and at least one $y$-region.
\item Let $Y(j) = 1$ if $R_j$ is a $y$-region and $Y(j) = 0$
otherwise.  Let $\CC_\sigma$ denote the relative cycle in $H_1 
(\nbd \cap \sing^{\critval + \ee} , \nbd \cap \sing^{\critval - \ee})$ 
represented by
$$\sum_{j=1}^k [ Y(j) - Y(j+1) ] \, \gamma_j$$
Then $\CC_\sigma$ represents the projection of $\alpha$ to 
this homology group, from which we easily recover $\paths (\sigma)$.
\item The cycle $\CC_* := \sum_{\sigma \in \Xi} \CC_\sigma$ represents
the projection of $\alpha$ to 
$H_1 (\sing^{\critval + \ee} , \sing^{\critval - \ee})$. 
\end{enumerate}
$\Cox$
\end{thm}

\noindent{\sc Proofs:} Let $c_1 > \cdots > c_\nu = c_{xy}$ be the
critical values of $h$ that are at least $c_{xy}$.  Let $\rind$ 
denote the set of $c \in \R$ for which $\partial \singx^{\geq c} 
\in \alpha$.  To establish~$(i)$ we will show:
\begin{enumerate}
\item[(a)] If $c \in \rind$ and $c_i > c > c_{i+1}$ then the entire
closed interval $[c_{i+1} , c_i ]$ is a subset of $\rind$;
\item[(b)] If $c \in \rind$ and $c > c_1$ then $[c_1 , \infty ) 
\subseteq \rind$;
\item[(c)] If $c_i \in \rind$ and $i < \nu$ then $c_i - \ee \in \rind$
for all sufficiently small $\ee > 0$;
\item[(d)] The set $\rind$ contains some $c > c_1$.
\end{enumerate}
By part~(i) of Theorem~\ref{th:morse}, the spaces 
$\{ \partial X^{\geq t} \; : \; cI > t > c_{i+1}$ 
are all  homotopic in $\singg$, hence represent the 
same homology class in $H_1 (\singg)$.  By continuity, 
this is true as well for $t = c_i$ and $t = c_{i+1}$, 
which proves~(a).  The proof of~(b) is identical.  

\begin{figure}[ht]
\centering
\includegraphics[scale=0.40]{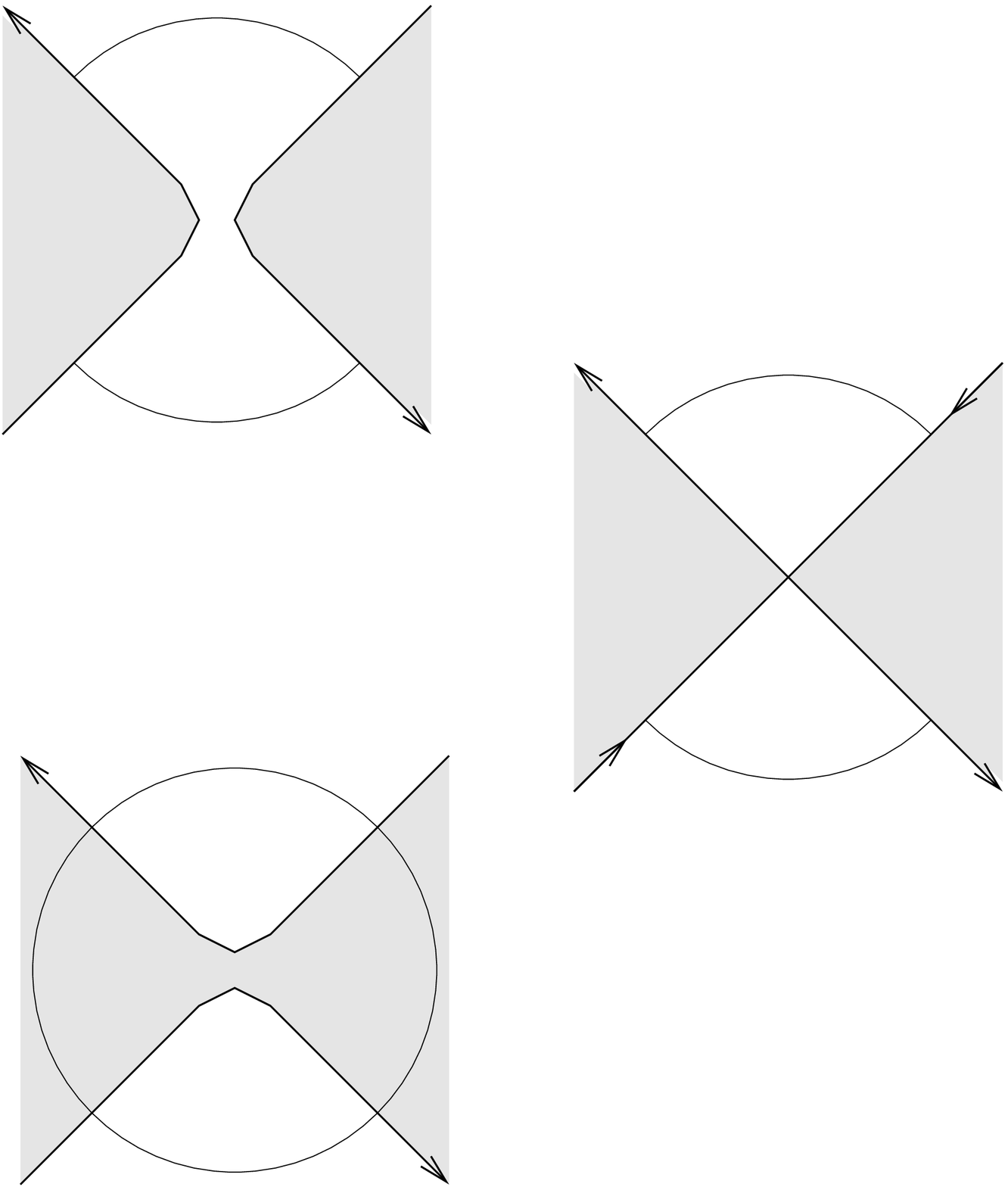}
\caption{$\sing^{\geq c}$ and its boundary for three values of $c$}
\label{fig:sectors}
\end{figure}

\begin{figure}[ht]
\centering
\includegraphics[scale=0.40]{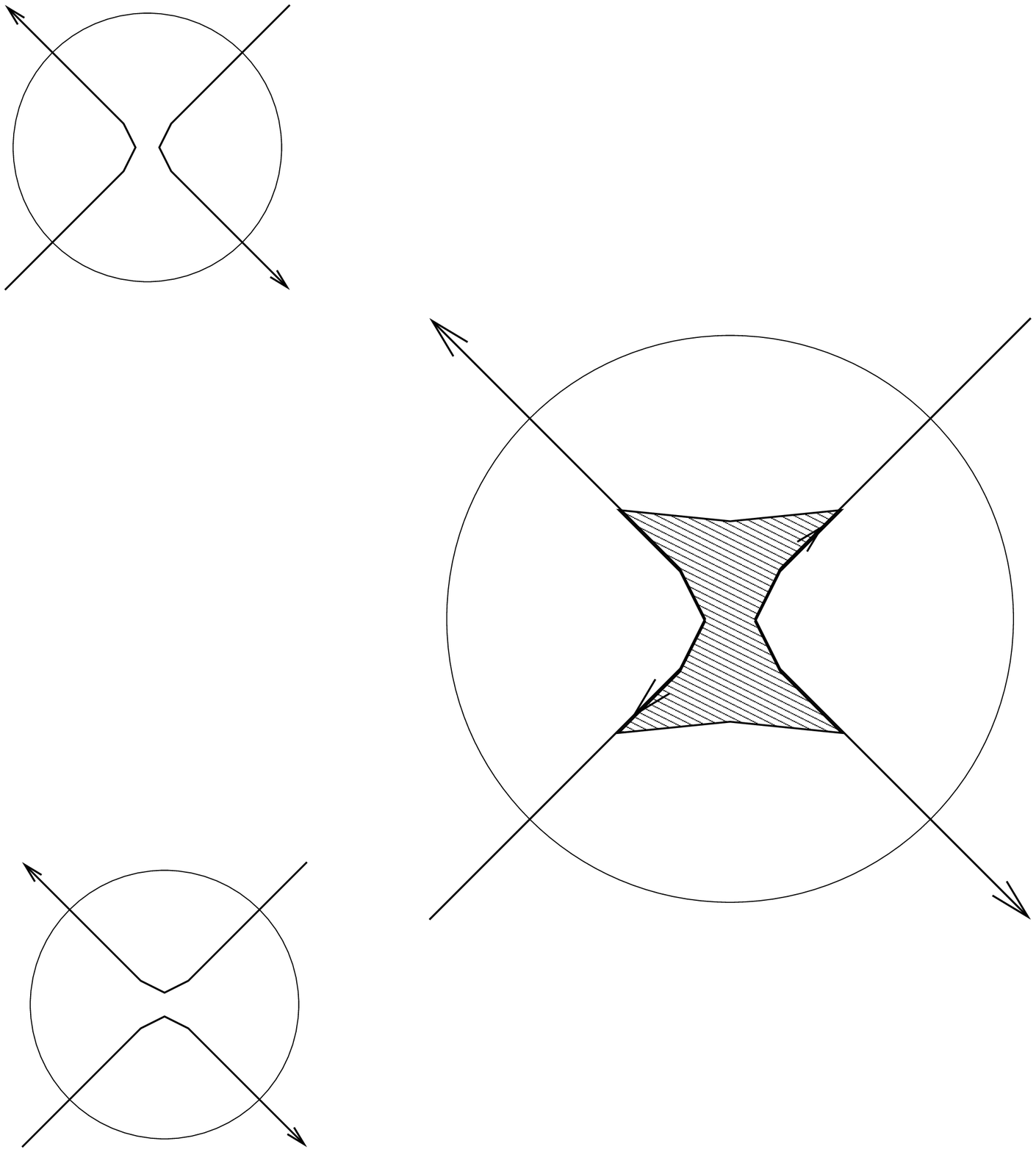}
\caption{$\sing^{\geq c_i + \ee}$ and $\sing^{\geq c_i - \ee}$ 
differ locally by a boundary}
\label{fig:polygon}
\end{figure}

To prove~(c), we first assume that there is a single critical
point $\sigma$ with $h(\sigma) = c_i$.  Let $k \geq 2$ be the
order to which the derivatives of $h$ vanish at $\sigma$.
Recall the neighborhood $U$ from Example~\ref{eg:harmonic}
and the analytic parametrization $\psi : U \to \sing$ with
$h(\psi (z)) = c_i + \Real \{ z^k \}$.  Locally, the 
image of $U$ is divided into $2k$ sectors with $h > c_i$ 
and $h < c_i$ in alternating sectors.  Figure~\ref{fig:sectors}
shows $\sing^{\geq c}$ (shaded) for three values of $c$
in the case $k=2$.  A circle is drawn to indicate a region
of parametrization for which $h = c_i + \Real \{ z^k \}$.
In the top diagram $c > c_i$, in the middle
diagram $c = c_i$ and in the bottom, $c < c_i$.  The arrows
show the orientation of $\partial \sing^{\geq c}$ inherited from
the complex structure of $\sing$.  The pictures for $k > 2$ are
similar but with more alternations.
Consider the first picture where $c > c_i$.  Because $c > c_{xy}$,
each of the $k$ shaded regions is in $\singx^{\geq c}$ or 
$\singy^{\geq c}$ but not both.  Let us term these regions 
``$x$-regions'' or ``$y$-regions'' accordingly.  Because $c_i > c_{xy}$, 
this persists in the limit as $c \downarrow c_i$, which means that
either all $k$ regions are in $\singy^{\geq c}$ or all $k$
regions are in $\singx^{\geq c}$.  In the former case, 
$\partial \singx^{\geq c}$ does not contain $\sigma$ for
$c$ in an interval around $c_i$ and the first Morse lemma
shows that $\partial \singx^{\geq c_i + \ee}$ is homotopic 
to $\partial \singx^{\geq c_i - \ee}$.  
In the latter case we consider the cycle 
$\partial \singx^{\geq c_i + \ee} + \partial B$ where 
$B$ is the polygon showed in Figure~\ref{fig:polygon}.
Because we added a boundary, this is homologous to 
$\partial \singx^{c_i + \ee}$.  But also it is homotopic to 
$\partial \singx^{c_i + \ee}$: within the parametrized neighborhood,
the lines may be shifted so as to coincide with $\partial 
\singx^{c_i + \ee}$, while outside this neighborhood the 
downward gradient flow provides a homotopy.  In the case where 
there is more than one critical point at height $c_i$, we may 
add the boundary of a polygon separately near each critical
$x$-point, that is, each point that is a limit as 
$\ee \downarrow 0$ of points in $\singx^{c_i + \ee}$.

Finally, we recall from Lemma~\ref{lem:transfer} that a cycle
$\CC$ representing $\alpha$ is constructed as $H_L [T_\ee] \cap \sing$.
Here $\ee$ is sufficiently small and $L$ sufficiently large that
this set is precisely $\{ (x,y) \in \sing : |x| = \ee \}$.  
By part~$(iv)$ of Theorem~\ref{th:classify}, the union of these
homotopic to $\partial \singx^{\geq c}$ for sufficiently large $c$.
Checking that the orientation given by Remark~\ref{rem:orient}
is the same as taking boundaries of regions oriented by the 
complex structure of $\sing$, we conclude~(d).  Part~$(i)$
of Theorem~\ref{th:critval} is a direct consequence of~(a)~--~(d).

Part~$(i)$ implies that for any $c \geq \critalt$, there is a 
cycle in the class $\alpha$ that is supported on $\sing^{\geq c}$, 
which shows that $\critval \leq \critalt$.  To see that $\critval 
\geq \critalt$, we suppose not and argue by contradiction.  The
construction in~\cite[Lemma~3.8]{pemantle-mvGF-AMS} shows that
for each $c_i > \critval$, there is a cycle representing $\alpha$
supported on $\sing^{\leq c_i + \ee}$ and any such cycle projects to
zero in the relative group $(\sing , \sing^{c_i - \ee})$.
For a contradiction, it suffices to check that 
$\partial \singx^{\critval + \ee}$ does not project to zero in 
$H_1 (\sing , \sing^{\leq \critval - \ee})$.  

Suppose first that there is a single critical point $\sigma$
at height $\critalt$.  Recall from Example~\ref{eg:harmonic}
that $H_1 (\sing^{\critalt + \ee} , \sing^{\leq \critalt - \ee})$
is generated by $\{ \gamma_i - \gamma_{i+1} 1 \leq i \leq k-1 \}$.
The only relation among these is that the sum of all of them
vanishes.  It follows that a sum of some subset of these 
vanishes if and only the subset is all or none.  Referring back
to Figure~\ref{fig:polygon}, it is clear that the cycle
$\partial \singx^{\geq \critalt}$ is represented in
$H_1 (\sing^{\critalt + \ee} , \sing^{\leq \critalt - \ee})$
by the sum of the $\gamma_{i+1} - \gamma_i$ over those $i$
for which the shaded sector between these two path segments
is an $x$-region.  
By the definition of $\critalt$, the point $\sigma$ is a limit
point of both $x$-regions and $y$-regions.  In other words,
the subset of shaded sectors that are in $\singx^{\geq \critalt}$
is a proper subset of all $k$ shaded sectors.  It follows that
$\partial \singx^{\geq \critalt}$ does not vanish in 
$H_1 (\sing , \sing^{\leq \critalt - \ee})$, and we have our contradiction.
This establishes~$(ii)$ of Theorem~\ref{th:critval} in the
case of only one critical point at height $\critval$
and also gives the characterization of $\CC_\sigma$ in
part~$(ii)$ of Theorem~\ref{th:Xi} in this case.

When there is more than one critical point $\sigma$ with 
$h(\sigma) = \critval$, it follows from the definition of
$\critalt$ that for at least one such $\sigma$, the $k$
local regions include at least one $x$-region and at least
one $y$-region.  The formula in part~$(ii)$ of Theorem~\ref{th:Xi}
continues to hold at each $\sigma$.  This vanishes when
all the regions are of one type, giving the characterization
of $\Xi$ in part~$(i)$ of Theorem~\ref{th:Xi}.  The direct
sum decomposition in part~$(ii)$ of Theorem~\ref{th:morse}
allows us to write $\CC_*$ as $\sum_{\sigma \in \Xi} \CC_\sigma$,
establishing part~$(iii)$.
$\Cox$

\section{Effective computation} \label{sec:computation}

By describing $\critval, \Xi$ and $\CC_*$ in terms of $x$- and 
$y$-regions, Theorems~\ref{th:critval} and~\ref{th:Xi} give us
a means to compute $\critval, \Xi$ and $\CC_*$.  One possible 
approach is as follows.  Instead of following $\partial \sing^{\geq c}$
down from $c = \infty$ to $\critval$, we follow paths upward
from each critical point until we find what we're looking for. 
A critical point of order $k$ has $k$ computable steepest 
ascent directions.  The backbone of our computation will 
therefore be:
\begin{enumerate}
\item Compute the critical points and order them by decreasing height.  
\item For the highest saddle, for each of the $k$ ascent paths,
check whether it converges to the $x$-axis or the $y$-axis.
\item if there is at least one of each type, set $\critval$ equal 
to the height of this saddle, set $\CC_\sigma$ according to 
part~$(ii)$ of Theorem~\ref{th:Xi}, and do the same for any 
other saddles at this height.  
\item If all are of the same type, then continue to the next lower 
saddle and iterate. 
\end{enumerate}
It should be clear that no new theorems are required in order
to implement this program.  However, a number of computational
apparati are required in order to make such a program completely 
effective; this section, which contains no major theorems, is
devoted to implementation.  

Chief among the needed computational apparati is a way to 
compute ascent paths.  Discretizing, we have the problem: 
given $(x,y) \in \sing$, 
produce a point $(x',y')$ with $h(x',y') > h(x,y)$ such that there 
is a strictly ascending path in $\sing$ from $(x,y)$ to $(x',y')$.  
Assuming we can do this, we must also ensure that this process ends 
after finitely many steps with identification as an $x$-component 
or a $y$-component.  This, together with some computer algebra, 
will complete our preparation for implementation.

\subsection{Ball arithmetic and {\sc Mathemagix}} \label{ss:ball}

For our implementation we have chosen the platform 
{\sc Mathemagix} \cite{vdH:mmx}.
This platform was designed for rigorous computations with 
objects of both algebraic and analytic nature.  On the one hand,
the algorithms in this paper indeed rely on purely symbolic 
pre-computations, which will be done using Gr\"obner bases.
On the other hand, our ascent paths are a special case of 
rigorous computations with analytic functions.  Interval arithmetic
is a systematic device for carrying out such
computations~\cite{Moor1966,AlHe1983,Neum1990,MKC2009,Rump2010}.
We will use a variant, called {\em ball arithmetic}~\cite{vdH2009}, 
which is more suitable for computations with complex numbers.
The use of interval arithmetic to achieve full rigor contrasts
to other homotopy continuation methods such as the package
Bertini~\cite{sommese-wampler} which use bootstrap
testing that is extremely reliable but still heuristic.
Although there is a large community for interval arithmetic, 
these techniques are not so common in other areas.
For this reason, we will recall the basic principles of ball arithmetic.

\subsubsection*{Ball arithmetic}

Given $c \in \C$ and $r \in \R^{\geq 0}$, denote by $\ball (c , r)$
the closed ball with center $c$ and radius $r$.  In what follows,
a complex \Em{ball number} is a ball $\ball (c , r)$ with $c \in 
\Q [i]$ and $r \in \Q^{\geq 0}$.  We think of such a ball number
as representing a generic complex number $z \in \ball (c,r)$.
The standard arithmetic operations can be defined on ball numbers
as follows:
\begin{eqnarray*}
\ball (c,r) + \ball (c',r') & = & \ball (c+c' , r+r') \\
\ball (c,r) - \ball (c',r') & = & \ball (c-c' , r+r') \\
\ball (c,r) \cdot \ball (c',r') & = & \ball (c \cdot c' , 
   (|c| + r) r' + r |c'|) 
\end{eqnarray*}
Given $* \in \{ + , - , \cdot \}$ and $\ball (c'', r'') = 
\ball (c,t) * \ball (c', r')$, these definitions have the
property that 
$$z \in \ball (c,r) \mbox{ and } z' \in \ball (c',r')
   \Longrightarrow z * z' \in \ball (c'',r'') \, .$$
Similar definitions can be given for other operations,
such as division, exponentiation, logarithm, and so forth.
Notice that any element $z \in \Q [i]$ gives rise to 
a ball number $\ball (z,0)$.

In practice, $\Q$ is usually replaced by $\float_{p,q}$,
the set of floating point numbers $x = m 2^e$ whose
mantissa $m \in \Z$ and exponent $e \in \Z$ have bounded 
bit lengths $p \in \N$ and $q \in \N$ respectively.
In that case the operations in $\float_{p,q}$ have 
additional rounding errors, and the definitions of the operations
on ball numbers must be adjusted to take into account these 
additional errors.  Furthermore, we must allow for the case 
where $r = +\infty$.  For more details, we refer 
to~\cite[Section~3.2]{vdH2009}.

Ball arithmetic allows for the reliable evaluation of more 
complicated expressions (or programs) $f(z_1 , \ldots , z_k)$
that are built up from the basic operations.  Indeed,
by induction over the size of the expression, it is easily
verified that
$$z_1 \in \ball (c_1 , r_1) , \ldots , z_k \in \ball (c_k , r_k)
   \Longrightarrow f(z_1 , \ldots , z_k) \in f (\ball (c_1 , r_1) ,
   \ldots , \ball (c_k , r_k)) \, .$$
Furthermore, it can be shown that ball arithmetic remains
continuous in the sense that the radius of $f (\ball (c_1 , r_1) ,
\ldots , \ball (c_k , r_k))$ will tend to zero if 
$r_1 , \ldots , r_k$ all tend to zero.

\subsubsection*{Complex algebraic numbers}

Thus far, we have no reliable zero test for ball numbers.  
In one direction, if $\ball (c,r) \cap \ball (c',r') = \emptyset$,
then we can be sure that $z \neq z'$ for any $z \in \ball (c,r)$
and $z' \in \ball (c',r')$. However the converse does not hold:
if $\ball (c,r) \cap \ball (c',r') \neq \emptyset$, then we do not
know whether $z = z'$ for $z \in \ball (c,r)$ and $z \in \ball
(c',r')$.  Of course, for some algorithms it is necessary to
be able to decide equality.  Fortunately, this can be accomplished
for the subclass of algebraic numbers.

Represent a complex algebraic number $z$ by a triple $(P,c,r)$
where $P \in \Q [t]$ is a square free polynomial and 
$\ball (c,r)$ is a ball number such that $z$ is the 
unique root of $P$ in $\ball (c,r)$.  This representation
is of course not unique.  The condition that $z$ be the
unique root of $P$ in $\ball (c,r)$ may be replaced by
a more explicit sufficient condition due to Krawczyk and 
Rump~\cite{Kraw1969,Rump1980,Rump2010}.  Consider the
expression
\begin{equation} \label{eq:KR}
\ball (c',r') = \Phi (\ball (c,r)) := 
   c - \frac{P(c)}{P'(c)} + \left ( 1 - \frac{P'(\ball (c,r))}{P'(c)}
   \right ) \cdot \ball (0,r),
\end{equation}
evaluated using ball arithmetic.  If $\ball (c',r')$ is contained
in the interior of $\ball (c,r)$ then it is certified that $P$ admits 
one and only one root $z \in \ball (c,r)$.  In particular, $z=0$ 
if and only if $P$ is divisible by $t$.  Moreover, for any $\ee > 0$,
there exists an iterate $\ball (c_n , r_n) = \Phi^n (\ball (c,r))$ 
such that $r_n < \ee$ and $z \in \ball (c_n , r_n)$.
In other words, the triple $(P,c,r)$ may be replaced by
a new triple $(P,c_n,r_n)$ for which the radius $r_n$ is 
arbitrarily small.

We will denote by $\algeb$ the set of complex algebraic numbers.
Using the above representation, arithmetic operations in $\algeb$
are easy to implement.  For instance, assume that we want to 
compute the sum of two algebraic numbers represented by 
$(P,c,r)$ and $(Q,c',r')$.  We first use symbolic algebra
to compute the square free part $R$ of the annihilator of
all sums of a root of $P$ and a root of $Q$.  Replacing
$\ball (c,r)$ and $\ball (c',r')$ by smaller balls if
necessary, we next ensure that $\ball (x'' , r'') :=
\ball (c,r) + \ball (c',r')$ satisfies the above
Krawczyk-Rump test for $R$.  The triple $(R,c'',r'')$
then represents the sum.

Combining Gr\"obner basis techniques, e.g.,~\cite{CLO2},
with algorithms for complex root finding of univariate
polynomials~\cite{Scho1982,BiFi2000}, the following can
be shown in a similar way: there exists an algorithm
which takes as input a zero-dimensional ideal of 
$\Q [ t_1 , \ldots , t_k]$ and produces as output 
all complex solutions in the form of vectors of
ball numbers describing disjoint poly-balls each
containing a unique root.

\subsubsection*{Ascent paths}

Given a point $z_0 \in \Q [i]$ with $f(z_0) \neq 0 \neq f' (z_0)$,
it will be important to be able to compute a line segment
of ascent, namely another rational complex number $z_1$
such that $|f|$ increases on the line segment from $z_0$
to $z_1$.  First, we note the following procedure to check
whether a ball lies entirely in the sector $S_{\pi/4}$ of
complex numbers whose arguments are strictly between $-\pi / 4$
and $\pi / 4$.  We consider a ball $\ball (c,r)$ that is known 
to intersect $S_{\pi/4}$ and we let $\ball (c' , r')$ denote
$\Real \{ \ball (c,r)^2 \}$ in ball arithmetic.  If $r' < |c'|$ 
then zero is not in $\ball (c' , r')$, hence by continuity
$\ball (c , r)^2$ lies in the open right half-plane.  
Because $\ball (c,r)$ is known to intersect $S_{\pi / 4}$,
it follows that $\ball (c,r)$ lies in $S_{\pi / 4}$ rather 
than $-S_{\pi/4}$.  In other words, to conclude that $\ball (c,r)
\subseteq S_{\pi/4}$, it is sufficient that $r' < |c'|$ in
the ball number that is the real part of the square of $\ball (c,r)$.

\begin{pr} \label{pr:ascent 0}
Let $f$ be locally analytic such that $f$ and $f'$ can be 
evaluated using ball arithmetic.  Let $z_0$ be a rational
point at which $f(z_0) \neq 0 \neq f' (z_0)$.  
\begin{enumerate}[(i)]
\item For $\ee > 0$ sufficiently small, the function
\begin{equation} \label{eq:g'}
u \mapsto \frac{f'(z_0 + u \frac{f(z_0)}{f'(z_0)} )}{f'(z_0)}
\end{equation}
evaluated at $u = \ball (0 , \ee)$ is contained in $S_{\pi / 4}$.
\item For such an $\ee > 0$, the function $|f|$ is strictly
increasing on the line segment from $z_0$ to $z_0 + \ee 
\frac{f(z_0)}{f'(z_0)}$.
\end{enumerate}
\end{pr}

\begin{remark} \label{rem:e}
Let $\ee_0 (f , z_0)$ denote the supremum of those $\ee > 0$
satisfying~\eqref{eq:g'} when evaluated in ball arithmetic.
Although implementation-dependent, this quantity is continuous
as a function of $z_0$. It is also easy to find an $\ee > 0$ which
satisfies~\eqref{eq:g'} and which is larger than $\ee_0 (f , z_0)/2$:
just keep doubling $\ee$ as long~\eqref{eq:g'} remains satisfied.
\end{remark}

\noindent{\sc Proof:} Change variables to the function $g$ defined by
$$g(u) := \frac{1}{f(z_0)} f \left ( z_0 + u \frac{f(z_0)}{f' (z_0)} 
   \right )  \, .$$
Thus $g(0) = g' (0) = 1$.  Because~1 is in the interior of 
$S_{\pi / 4}$, continuity of ball arithmetic implies that 
$g' (\ball (0 , \ee))$ lies in $S_{\pi / 4}$ for sufficiently
small positive $\ee$.  The function in~\eqref{eq:g'} is
precisely $g'$, which proves~$(i)$.  For~$(ii)$, observe that
for $u \in [0,\ee]$ we have
\begin{equation} \label{eq:g}
g(u) = \int_0^u g' (v) \, dv \in S_{\pi/4} \, .
\end{equation}
This implies $\Arg (g'(u) / g(u)) \in (-\pi/2 , \pi/2)$
for all $u \in [0,\ee]$.  Thus $|g|$ increases strictly 
on $[0,\ee]$, and changing variables back to $f$ gives~$(ii)$.
$\Cox$

More generally, we will need to compute a segment of ascent
for algebraic points where the derivative vanishes to some 
order $k$.  The same argument as for Proposition~\ref{pr:ascent 0}
easily shows the following.
\begin{pr} \label{pr:saddle ascent 0}
Let $f$ be locally analytic such that $f$ and its first $k$
derivatives may be evaluated using ball arithmetic.  Let $z_0$ 
be a rational point at which $f(z_0) \neq 0 \neq f^{(k)} (z_0)$
while $f^{(1)} (z_0) = \cdots = f^{(k-1)} (z_0) = 0$.  Define
$$g(u) := \frac{1}{f(z_0)} f \left ( z_0 + u \left [ \frac{f(z_0)}{f' (z_0)} 
   \right ]^{1/k} \right ) $$
where any choice of the $k^{th}$ root is allowed.  Then when 
$\ee > 0$  is sufficiently small, $g^{(k)} (\ball (0 , \ee))
\subseteq S_{\pi/4}$ and for such an $\ee$, the magnitude
of $f$ will increase strictly on the line segment from $z_0$
to $z_0 + \ee [f(z_0) / f^{(k)} (z_0)]^{1/k}$.  Furthermore,
computing such an $\ee$ for each choice of $1/k$ power and 
taking the minimum ensures that $|f|$ increases on all $k$ 
line segments simultaneously.
$\Cox$
\end{pr}

\subsection{The use of ascent paths to compute invariants on $\sing$}

We now harness these computational devices to compute
the topological invariants of $\sing$ that are required
for asymptotics, namely $c_*, \Xi$ and $\paths$.  We begin
by transferring ascent segments of a function 
(Propositions~\ref{pr:ascent 0} and~\ref{pr:saddle ascent 0})
to ascending paths on the Riemann surface $\sing$.
Let $\xns$ denote the set of points $\zz = (x , y) \in 
\sing \setminus \saddles$ with $x \in \Q [i]$ and
$\partial Q / \partial x (\zz) \neq 0$.  We may represent
elements of $\xns$ as ball numbers with first coordinate
radius zero and second coordinate radius small enough that
the ball contains only one root of $Q(x,\cdot)$. 
\begin{pr}[rigorous ascent step] \label{pr:ascent}
There is a ball-computable function $\phi : \xns \to \xns$
with the following properties.  Let $\zz_0 = (x_0 , y_0) \in \xns$
and denote $\zz_1 := (x_1 , y_1) := \phi (\zz_0)$.
Then the line segment $[x_0 , x_1]$ lifts uniquely to 
a curve in $\sing$ connecting $\zz_0$ to $\zz_1$, along 
which $h$ is strictly increasing.  Furthermore, $\phi$
may be chosen so that $h(\zz_1) - h(\zz_0)$ is bounded
below by a positive constant $c_K$ on any bounded subset $K$ 
of $\xns$ whose closure avoids $\saddles$ and points where
$\partial Q / \partial x$ vanishes.
\end{pr}

\noindent{\sc Proof:} Given $\zz_0 \in \xns$, let $y$ be
the locally analytic function such that $y(x_0) = y_0$
and $Q(x,y(x)) = 0$.  Apply Proposition~\ref{pr:ascent 0}
to $f (z) := \exp [h (z , y(z))]$ and $z_0 = x_0$, obtaining
a segment $[x_0 , x_1] = [z_0 , z_0 + \ee_1 f(z_0) / f'(z_0)$.
It is easy, computationally, to choose $\ee_1$ always to be
at least $\ee_0 / 2$.  Define $\phi_1 (\zz_0) := \zz_1 := 
(x_1 , y_1)$.  The lifting of $[x_0,x_1]$ to $\sing$ is a 
path along which $h$ increases, hence $h(\phi_1 (\zz)) 
- h(\zz)$ is strictly positive.  Let us check that this 
difference is bounded below by a positive constant on 
compact sets.  

We know $\ee_1$ is bounded away from zero
on compact sets avoided by $\saddles$ and $\partial Q / 
\partial x$ because $\ee_0$ is.  Letting $h (z_0 , t)$ 
denote the derivative of $|f|$ along the lifted segment at
$z_0 + t f(z_0) / f'(z_0)$, and observing that $h(z_0,0)$
is bounded below by a positive constant on compact
sets avoiding $\saddles$, we conclude that 
$h(\phi_1 (\zz)) - h(\zz)$ is indeed
bounded below by a positive constant $c_K'$
on compact sets, $K$.  

This gives the conclusion we desire except that $x_1$ 
may not be rational.  To correct this, note that 
all we used about $\ee_1$ was that was guaranteed
to be in the interval $[c \ee_0 , \ee_0]$ for some
constant $c$.  Having computed $\ee_1$, we may
choose $\ee > 0.99 \ee_1$ such that $\ee f(z_0) / f'(z_0)$
is rational.  Denoting the endpoint of the new shorter
lifted path by $\phi (\zz)$, the same argument as before
now shows that $h(\phi (\zz)) - h(\zz)$ is bounded
from below by a positive constant, finishing the proof.
$\Cox$

Having defined a step from a point not near a saddle, we
next find a way to ascend out of a saddle.  The following
proposition describes a way to do this.  The proof is
nearly identical to the proof of Proposition~\ref{pr:ascent}
and is omitted.
\begin{pr}[ascent from a saddle] \label{pr:saddle ascent}
Let $\zz_0 = (x_0 , y_0) \in \saddles$ with $x_0 \in \Q [i]$.
Suppose $\partial Q / \partial x$ does not vanish at $\zz_0$.  
Then we may compute a rational $\psi (x_0)$ such that
the union of the radial line segments 
$$\left \{ [x_0 , x_0 + \phi (x_0) e^{2 \pi i j / k}] : 0 \leq j \leq k-1 
   \right \} $$
lifts uniquely to a union of paths from $\zz_0$ on $\sing$ 
on each of which $h$ is strictly increasing.
$\Cox$
\end{pr}

The following proposition allows us to terminate an ascent path
when it comes close enough to the $x$-axis or $y$-axis.
\begin{pr}[arrival at $\singx$ or $\singy$] \label{pr:pole nbd}
Fix $\rhat$ with $r/s \notin \bad$ and let $h = h_{\rhat}$.  
Let $\cmax$ denote the greatest critical value of $h$ and
$\cmin$ denote the least critical value.
Let $\ee$ be small enough so that there are no critical 
points $(x,y)$ with $|x| < \ee$.  Then any point $(x,y) 
\in \sing$ with $|x| < \ee$ and $h(x,y) > \cmin$ is in
$\singx^{> \cmax}$.
\end{pr}

\begin{unremark}
The hypothesis is formulated so as to be easily checked using
ball arithmetic: letting $J$ be the ideal generated by the
critical point equations~\eqref{eq:crit}, we take $\ee$ 
to be sufficiently small such that the ball solutions to $J$ 
are disjoint from the ball $\ball (0,\ee)$.  Also note that the roles 
of $x$ and $y$ may be switched to yield an analogous result.
\end{unremark}

\noindent{\sc Proof:} Let $\ee$ be as in the hypothesis.  We know 
that the graph of $\sing$ over $\ball (0,\ee)$ is a union of graphs of
Puiseux series.  Recall that as $x \to 0$ on one of these components,
$h \to \infty$ except possibly for some ``low'' components on which 
$y$ is unbounded and $h \to -\infty$.  By hypothesis, $h$ cannot take 
on a critical value on any of these components, so we have $h > \cmax$
on the high components and $h < \cmin$ on the low components.
This accounts for all solutions to $|x| < r$ on $\singg$, hence
$|x| < r$ together with $c > \cmin$ suffices to assure that
a point $(x,y) \in \sing$ is in $\singx^{> \cmax}$.
$\Cox$

The next proposition allows ascent paths to terminate when
they come sufficiently near a saddle.  Without this 
improvement there is a danger that an ascent path could
converge to a saddle along an infinite sequence of ever 
smaller steps.  This necessary improvement is also a big
time-saver.  Step~2 of the algorithm at the beginning of
Section~\ref{sec:computation} calls for us to identify
which axis is the limit of each ascent path from a given
saddle.  Proposition~\ref{pr:saddle nbd} classifies an
ascent path once it comes near any higher saddle, thereby
saving the remainder of the journey to the axes.
\begin{pr}[bypassing a saddle] \label{pr:saddle nbd}
Let $c > \critval$ be a critical value of $h$ and let 
$(x_0,y_0)$ be a critical point at height $c$ in 
$\singx^{\geq c}$.  Suppose there is an $\ee > 0$ 
for which the following conditions hold.
\begin{enumerate}[(i)]
\item For each $x$ such that $|x - x_0| \leq \ee$ there is at
most one solution to $Q(x,y) = 0$ with $|y - y_0| \leq \ee$;
\item $(x_0 , y_0)$ is the unique critical point $(x,y)$ of $h$ 
on $\sing$ with $|x - x_0|  , |y - y_0| \leq \ee$.
\end{enumerate}
Then $|x - x_0| , |y - y_0| < \ee$ and $Q(x,y) = 0$ together 
imply $(x,y) \in \singx^{> \critval}$.  
\end{pr} 

\begin{unremarks}
An identical result holds if $X^{\geq c}$ and $X^{> c_*}$ are 
replaced by $Y^{\geq c}$ and $Y^{> c_*}$ throughout.  Also,
hypotheses~$(i)$ and~$(ii)$ are easily checked in ball arithmetic.
\end{unremarks}

\noindent{\sc Proof:} Our first hypothesis on $\ee$ implies 
that the set 
$$S := \sing \cap \{ |x - x_0| < \ee \} \cap \{ |y - y_0| < \ee \}$$
is a graph over $\ball (x_0 , \ee)$ of a univalent holomorphic function.
Our second hypothesis implies that $\sing^{= \critval}$ does
not intersect $S$.  Because $\sing^{= \critval}$ is the boundary 
of $\sing^{< \critval}$ this implies that the entire set $S$ is
contained in one component of $\sing^{> \critval}$.  Because
$(x_0,y_0)$ is an element of $S$ and is also in $\singx^{> \critval}$,
we conclude that $S \subseteq \singx^{> \critval}$ which is the
desired conclusion.  
$\Cox$

Let $\sigma_1 , \sigma_2 , \ldots, \sigma_m$ be ball algebraic 
representations of the critical points in weakly descending 
order of height.  Fix $\rhat \notin \bad$.  The algorithm
outlined at the beginning of Section~\ref{sec:computation}
may now be described in detail, proved to terminate and
proved to produce a correct answer.  In the end, the index
$j_0$ will be such that $c_{j_0} = \critval$; each critical
point will be classified as $x$, $y$, mixed or unclassified; 
the set $\Xi$ will be the set of mixed saddles; the unclassified
saddles will be the ones lower than $\critval$.
\begin{alg}[Computation of $\critval, \Xi$ and $\CC_*$]
\label{alg:1} ~~\\[-5ex]
\begin{enumerate}
\item {\bf (initialization)} 
\begin{itemize}
\item Set $j := 1$.  
\item Let $\ee$ be as in Proposition~\ref{pr:pole nbd} and let
$\ee'$ be the analogous quantity with $x$ and $y$ switched.
\item For $1 \leq j \leq m$ set $\class (\sigma_j) := $
``unclassified''. 
\item Set $\critval := -\infty$.
\end{itemize}
\item {\bf (main loop)} Repeat while $h(\sigma_j) \geq \critval$:
\begin{itemize}
\item Set $\zz_0 := \sigma_j$.
\item Apply Proposition~\ref{pr:saddle ascent} with 
$(x_0 , y_0) = \zz_0$ for each of the $k$ choices of $1/k$ power 
and let $\zz_i$ for $1 \leq i \leq k$ denote the 
other endpoint of the resulting line segment.
\item For $i=1$ to $k$, initialize $\classify_i$ to ``unclassified''.
\item For $i=1$ to $k$, initialize $\done$ to False and 
repeat until $\done$:
\begin{itemize}
\item If the $x$-component of $\zz_i$ has modulus less than 
$\ee$ then set $\classify_i := $``$x$'' and $\done := $ True;
\item else if the $y$-component of $\zz_i$ has modulus less than 
$\ee'$ then set $\classify_i := $``$y$'' and $\done := $ True;
\item else if $\zz_i$ is in the ball $\ball (\sigma_p)$ in the conclusion of 
Proposition~\ref{pr:saddle nbd} for any $(x_0 , y_0) = \sigma_p : 
p < j$ then set $\classify_i := \classify_p$ and $\done := $ True;
\item else set $\zz_i := \phi (\zz_i)$ as in Proposition~\ref{pr:ascent}.
\end{itemize}
\item 
\begin{itemize}
\item If $\class_i = $``$x$'' for all $1 \leq i \leq k$ then set
$\class (\sigma_j) = $``$x$''.
\item else if $\class_i = $``$y$'' for all $1 \leq i \leq k$ then set
$\class (\sigma_j) = $``$y$''.
\item else set $\class (\sigma_j) := $ ``mixed''.
\end{itemize}
\item If $\critval = $ ``undefined'' and $\class (\sigma_j) = $ ``mixed''
then set $\critval := h (\sigma_j)$.
\item $j := j+1$.
\end{itemize}
\item {\bf (finishing the computation)}
\begin{itemize}
\item Set $\Xi = \emptyset$.
\item For $i = 1$ to $j$, if $\class (\sigma_i) = $``mixed'' then do: 
\begin{itemize}
\item $\Xi := \Xi \cup \{ \sigma_j \}$.
\item Compute $\CC_{\sigma_j}$ via Theorem~\ref{th:Xi} with
$Y(i) = \one_{\classify_i = y}$.
\end{itemize}
\end{itemize}
\end{enumerate}
\end{alg} \vspace{0.12in}

\begin{thm} \label{th:path}
This algorithm terminates and correctly computes $\critval,
\Xi$ and $\{ \C_\sigma : \sigma \in \Xi \}$.
\end{thm}

\noindent{\sc Proof:} The sequence $\zz_0 , \zz_i , \phi (\zz_i) , 
\phi (\phi (\zz_i)), \ldots$ defines a polygonal path in $x$
whose lifting to $\sing$ is an ascent path for $h$ along 
which  the value of $h$ on each step increases by an amount 
bounded away from zero.  We must therefore reach the
terminal conditions of the $\done$ loop in a number of steps
bounded by $L$, where $L = M - \cmin$ divided by this minimum
step increase and $M$ is the greatest height of any $(x,y) \in 
\sing$ with $|x| = \ee$.  The $k(\sigma)$ ascent paths from each 
critical point $\sigma$ verify whether each of the $k$ regions
is and $x$-region or a $y$-region; the algorithm computes 
$\critalt$ and thus, by Theorem~\ref{th:critval}, computes $\critval$.
Theorem~\ref{th:Xi} shows that the remainder of the algorithm
computes $\Xi$ and $\CC_*$.
$\Cox$

\subsection{Uniform asymptotics as $\rhat$ varies}

Let 
$$\degen := \{ \rhat : h_{\rhat} 
   \mbox{ has a degenerate critical point } \} \, $$
where degenerate means having order~3 or higher.  
We use this name because the degenerate saddle of order~3
is called a ``monkey saddle'' (downward regions for
both legs and the tail).
\begin{pr} \label{pr:monkey}
Assume for nondegeneracy that $Q$ is not a binomial.
Let $J$ be the ideal in $\C [x , y , \lambda]$ generated
by the following three polynomials. 
\begin{eqnarray}
&& Q \nonumber \\
&& x Q_x - \lambda y Q_y \nonumber \\
&& y Q_y^2 (Q_x + x Q_{xx}) + x Q_x^2 (Q_y + y Q_{yy}) 
   - 2xy Q_x Q_y Q_{xy} \label{eq:psi}
\end{eqnarray}
Then $(r,s) \in \degen$ implies $r/s$ is the $\lambda$-coordinate 
of a solution to $J$.  In particular, the set $\degen$ is the 
intersection of $(0,\infty)$ with the zero set of the
elimination polynomial $f$ of $\lambda$ in $J$.
\end{pr}

\noindent{\sc Proof:} The first two generating polynomials
for $J$ define an algebraic function 
\begin{equation} \label{eq:lambda}
\lambda (x,y) := \frac{x Q_x}{y Q_y}
\end{equation}
on the variety $\{ Q = 0 \}$.  This is the \Em{direction} function
of~\cite{PW1} in the sense that a point $(x,y) \in \sing$ is 
critical for $h_{\rhat}$ if and only if $\lambda (x,y) = r/s$.
The derivative of $\lambda$ on $\sing$ with respect to, say, $x$,
is given by
$$\left. \frac{d}{dx} \right |_\sing \; (\lambda) = \lambda_x 
   + \lambda_y \frac{-Q_x}{Q_y}$$
which vanishes exactly when the Wronskian $\lambda_x Q_y 
- \lambda_y Q_x$ vanishes.  Using the explicit formula~\eqref{eq:lambda}
for $\lambda$ shows the Wronskian, up to some factors of $y$ and $Q_y$,
to be the third of the given generators of $J$.  Vanishing of the
derivative of $\lambda$ on $\sing$ is a necessary condition for
the coalescing of two solutions to $\lambda = c$ on $\sing$.  This
proves the main conclusion.  If there is no elimination polynomial
for $f$ in $J$ then $\lambda$ is constant on $\sing$ which implies
that $Q$ is binomial.  We have assumed not, from which the last 
statement follows.  
$\Cox$

Let $\At := \bad \cup \degen$ denote the set of troublesome 
directions.  We see from Proposition~\ref{pr:monkey} and the
definition of $\bad$ that $\At$ is readily computed.  Its 
complement in the arc of $\rhat$ in the positive quandrant 
is a finite union of intervals with endpoints whose slopes 
are algebraic.  On each interval in $\At^c$, the asymptotic 
estimate for $a_\rr$ in Lemma~\ref{lem:integrals} is uniform 
over compact sub-intervals.  This allows us effectively to 
give uniform asymptotics in all nondegenerate directions.
\begin{alg} \label{alg:2} ~~\\[-5ex]
\begin{itemize}
\item Compute $\At$
\item Enumerate the intervals of $\At^c$ in the open arc 
parametrized by $\lambda \in (0,\infty)$ and choose one 
rational value $r(I)$ in each interval $I$.
\item for each interval $I$ do:
\begin{itemize}
\item Set $\rhat = r(I)$
\item Use Algorithm~\ref{alg:1} to compute $c_* , \Xi$ and
$\CC_*$ 
\item Use Lemma~\ref{lem:integrals} to write a uniform estimate
for $a_\rr$ with $\rhat \in I$ by summing~\eqref{eq:integrals} over
$\sigma \in \Xi$ and $i \in \paths (\sigma)$.
\end{itemize}
\item Interpret the formula as holding uniformly over compact
subintervals of $I$, where multi-valued quantities depending
on $\rhat$ are extended by homotopy from their values at $r(I)$.
\end{itemize}
\end{alg}

\subsubsection*{Discussion of remaining directions}

To complete the asymptotic analysis of $a_\rr$, we need to 
know what happens in the remaining directions.  For each 
$\rhat \in \degen$, we have a formula given by inserting 
the estimates from Lemma~\ref{lem:integrals} into 
Algorithm~\ref{alg:1}.  These estimates are known to 
extend to Airy functions on the rescaled window
$r = \lambda s + O(s^{1/2})$.  Because the focus of this
paper is how to make the computation of $\Xi$ and $\paths$
effective, we do not discuss the extensions of~\eqref{eq:integrals}
to Airy function estimates in the present work.  

When $r/s  = \lambda_0 \in \bad$, 
there is a projective point of $\sing$ of finite height.  It is possible
that this point is a smooth point and not a critical point, in which
case a change of chart maps gets rid of it; the formulae for the
two adjacent intervals $I$ and $I'$ will agree and will be valid throughout
the union $I \cup I' \cup \{ \lambda_0 \}$.  More often, however,
the projective point fails to be a smooth point.  In this case
the formulae on $I$ and $I'$ will in general be different and 
asymptotics in the direction $\lambda_0$ will be given by something
other than~\eqref{eq:integrals}.  

\bibliographystyle{alpha}
\bibliography{RP}

\end{document}